\documentstyle[12pt]{article}

\def\no{\nonumber}
\def\be{\begin{equation}}
\def\ee{\end{equation}}
\def\ba{\begin{eqnarray}}
\def\ea{\end{eqnarray}}
\def\tilde{\widetilde}
\def\btu{\Delta}
\def\btd{\nabla}
\def\e1{\epsilon}
\def\o1{\omega}
\def\01{\Omega}
\def\c1{\gamma}
\def\al{\alpha}

\def\g1{\Sigma}
\def\lmd{\lambda}
\def\l1{\Lambda}
\def\v1{\varphi}
\def\d1{\delta}
\def\part{\partial}
\def\fr{\frac}
\def\f2{F}
\def\h2{{\bf H}}
\def\a2{{\bf A}}
\def\x2{{\bf X}}
\def\t1{\theta}
\def\b1{\beta}

\def\bs{\begin{eqnarray*}}
\def\es{\end{eqnarray*}}
\def\p{\partial}
\def\m1{\Theta}
\def\w1{\wedge}
\def\la{\langle}
\def\ra{\rangle}

\def\O{\Omega}
\def\o{\omega}
\def\S{\Sigma}

\begin{document}

\title{Singularity of Mean Curvature Flow of
Lagrangian Submanifolds}

\author{\begin{tabular}{ccc}
Jingyi Chen\thanks{Chen is supported partially by a Sloan fellowship and a NSERC grant.} &Jiayu
Li\thanks{Li is partially supported by the National Science Foundation of China
and by the Partner Group of MPI for Mathematics.}\\
Department of Mathematics & Institute of Mathematics\\
University of British Columbia & Fudan University, Shanghai\\
Vancouver, B.C. V6T 1Z2 & Academia Sinica, Beijing\\
Canada & P.R. China\\
jychen@math.ubc.ca&lijia@math.ac.cn
\end{tabular}}

\date{}

\maketitle

\abstract{In this article we study the tangent cones at
first time singularity of a Lagrangian mean curvature flow.
 If the initial compact submanifold $\S_0$
is Lagrangian and almost calibrated by ${\mbox{Re}}\,\O$ in
a Calabi-Yau $n$-fold $(M,\O)$, and
$T>0$ is the first blow-up time of the mean curvature flow, then the
tangent cone of the mean curvature flow at a singular point $(X_0,T)$ is a
stationary Lagrangian integer multiplicity current in ${\bf R}^{2n}$ with volume
density greater than one at $X_0$. When $n=2$, the tangent cone
consists of a finite union of more than one 2-planes in ${\bf
R}^4$ which are complex in a complex structure on ${\bf R}^4$.
}

\newtheorem{theorem}{Theorem}[section]
\newtheorem{lemma}[theorem]{Lemma}
\newtheorem{corollary}[theorem]{Corollary}
\newtheorem{remark}[theorem]{Remark}
\newtheorem{definition}[theorem]{Definition}
\newtheorem{proposition}[theorem]{Proposition}
\newtheorem{example}[theorem]{Example}

\section{Introduction}

Let $M$ be a compact Calabi-Yau manifold of complex dimension $n$
with a K\"ahler form $\o$, a complex structure $J$, a K\"ahler
metric $g$ and a parallel holomorphic $(n,0)$-form $\O$ of unit
length. An immersed submanifold $\S$ in $M$ is Lagrangian if
$\o|_\S=0$. The induced volume form $d\mu_\S$ on a Lagrangian
submanifold $\S$ from the Ricci-flat metric $g$ is related to $\O$
by
\begin{equation}
\O|_\S= e^{i\theta}d\mu_\S =\cos\theta d\mu_\S+i\sin\theta d\mu_\S,
\end{equation}
where the phase function $\theta$ is multi-valued and is well-defined up to an
additive constant $2k\pi,k\in {\mathbf Z}$. Nevertheless, $\cos\theta$ and $\sin\theta$
are single valued function on $\S$. For any tangent vector $X$ to $M$
a straightforward calculation shows
\begin{equation}
X\theta = -g(\h2,JX)
\end{equation}
where $\h2$ is the mean curvature vector of $\S$ in $M$ (cf. [HL],
[TY]). Equivalently, $\h2=J\nabla\theta$. The Lagrangian
submanifold $\S$ is {\it special}, i.e. it is a minimal
submanifold, if and only if $\theta$ is constant. When $\theta$ is
constant on a Lagrangian submanifold $\S$, the real part of
$e^{-i\theta}\O$ is a calibration of $M$ with comass one and $\S$
is a volume minimizer in its homology class [HL]. Let ${\mbox
Re}\O$ be the real part of $\O$. A Lagrangian submanifold is
called {\it almost calibrated} by ${\mbox Re}\O$ if $\cos\theta >
0$.

Constructing minimal Lagrangian submanifolds is an important
but very challenging task. In a compact K\"ahler-Einstein surface,
Schoen and Wolfson [ScW] have shown the existence of a branched surface which minimizes
area among Lagrangian competitors in each Lagrangian homology class, by variational method.

For a one-parameter family of immersions $\f2_t=\f2(\cdot ,t):\g1\to M$, we denote
the image submanifolds by $\g1_t=\f2_t(\g1)$. If $\S_t$ evolves along the gradient flow
of the volume functional, the first variation of the volume functional asserts that
$\g1_t$ satisfy a mean curvature flow equation:
\begin{equation}\label{meaneqn}
\left\{\begin{array}{lll}
\displaystyle \fr{d}{dt}\f2(x,t)=\h2(x,t)\\
\displaystyle \,\,\,\,\,\,\f2(x,0)=\f2_0(x),
\end{array}\right.
\end{equation}

When $\S$ is compact the mean curvature flow (\ref{meaneqn})
has a smooth solution for short time $[0,T)$ by the standard parabolic theory. If $\S_0$ is
Lagrangian in a K\"ahler-Einstein ambient space $M$, Smoczyk has shown that $\S_t$
remains Lagrangian for $t<T$ and the phase function $\theta$ evolves by
\begin{equation}\label{theta}
\frac{d\theta}{dt}=\Delta\theta
\end{equation}
where $\Delta$ is the Laplacian of the induced metric on $\S_t$
([Sm1-3], also see [TY] for a derivation of (\ref{theta})). It
then follows that
\begin{equation}\label{cos}
\frac{\partial\cos\theta}{\partial t}=\Delta\cos\theta + |\h2|^2\cos\theta.
\end{equation}
If the initial Lagrangian submanifold $\S_0$ is almost calibrated, $\S_t$ is almost calibrated,
i.e. $\cos\theta > 0$,
along a smooth mean curvature flow by the parabolic maximum principle.

It is well-known that if $|\a2|^2$, where $\a2$ is the second fundamental
form on $\g1_t$, is bounded uniformly as $t\to T>0$ then (\ref{meaneqn}) admits a
smooth solution over $[0,T+\e1)$ for some $\e1>0$. When
$\max_{\g1_t} |\a2|^2$ becomes unbounded as $t\to T$,
we say that the mean curvature flow develops a singularity at $T$. A lot of work
has been devoted to understand these singularities (cf. [CL1-2],
[E1-2], [H1-3], [HS1-2], [I1], [Wa], [Wh1-3].)

In this paper, we shall study the tangent cones at singularities of the
mean curvature flow of a compact Lagrangian submanifold in a compact
Calabi-Yau manifold. Especially, we shall focus on the
structure of tangent cones of the mean curvature flow where a
singularity occurs at the first singular time $T<\infty$.

To describe the tangent cones, suppose that $(X_0,T)$ is a
singular point of the flow (\ref{meaneqn}), i.e. $|\a2(x,t)|$
becomes unbounded when $(x,t)\to (X_0,T)$. For an arbitrary
sequence of numbers $\lmd\to\infty$ and any $t<0$, if
$T+\lambda^{-2}t>0$ we set
$$
F_\lmd(x,t)=\lmd (F(x,T+\lmd^{-2}t)-X_0).
$$
We denote the scaled submanifold by $(\g1_t^\lmd,d\mu_t^\lmd)$. If
the initial submanifold is Lagrangian and almost calibrated by
${\mbox{Re}}\,\O$, it is proved in Proposition \ref{exit} that
there is a subsequence $\lmd_i\to\infty$ such that for any $t<0$,
$(\g1_t^{\lmd_i},d\mu_t^{\lmd_i})$ converges to
$(\g1^{\infty},d\mu^{\infty})$ in the sense of measures; the limit
$\g1^\infty$ is called a {\it tangent cone arising from the
rescaling $\lmd$}, or simply a {\it $\lmd$ tangent cone at
$(X_0,T)$}. This tangent cone is independent of $t$ as shown in
Proposition \ref{exit}.

There is also a time dependent scaling which we would like to consider
\begin{equation}
\tilde{F}(\cdot ,s)=\fr{1}{\sqrt{2(T-t)}}F(\cdot ,t),
\end{equation}
where $s=-\fr{1}{2}\log (T-t)$, $c_0\leq s<\infty $. Here we have
chosen the coordinates so that $X_0=0$. Rescaling of this type
arises naturally in classification of singularities of mean
curvature flows [H2]: assume $ \lim_{t\to
T^-}\max_{\g1_t}|\a2|^2=\infty$, if there exists a positive
constant $C$ such that $ \limsup_{t\to T^-}\left((T-t)\max_{\g1_t}
|\a2|^2\right)\leq C, $ the mean curvature flow $\f2$ has a Type I
singularity at $T$; otherwise it has a Type II singularity at $T$.
Denote $\tilde{\g1}_s$ the rescaled submanifold by
$\tilde{F}(\cdot,s)$. If a subsequence of $\tilde{\g1}_s$
converges in measures to a limit $\tilde{\g1}_\infty$, then the
limit is called a {\it tangent cone arising from the time
dependent scaling at $(X_0,T)$}, or simply a {\it $t$ tangent
cone}. In this paper, a {\it tangent cone}  of the mean curvature
flow at $(X_0,T)$ means either a $\lmd$ tangent cone or a $t$
tangent cone at $(X_0,T)$.

The main result of this paper is

\begin{theorem}\label{main} Let $(M,\O)$ be a compact Calabi-Yau manifold of
complex dimension $n$. If the initial compact submanifold $\S_0$
is Lagrangian and almost calibrated by ${\mbox{Re}}\,\O$, and
$T>0$ is the first blow-up time of the mean curvature flow
(\ref{meaneqn}), and $(X_0,T)$ is a singular point, then the
tangent cone of the mean curvature flow at $(X_0,T)$ is a
stationary Lagrangian integer multiplicity current in ${\bf
R}^{2n}$ with volume density greater than one at $X_0$. When
$n=2$, the tangent cone consists of a finite union of more than
one 2-planes in ${\bf R}^4$ which are complex in a complex
structure on ${\bf R}^4$.
\end{theorem}

For symplectic mean curvature flow in K\"ahler-Einstein surfaces, results similar to
Theorem \ref{main} was obtained in [CL1].
The authors are grateful to Professor Gang Tian for stimulating conversations.

\section{Existence of $\lmd$ tangent cones}

This section contains basic formulas and estimates which are
essential for this article. First, we will derive a monotonicity
formula which has a weight function introduced by the $n$-form
${\mbox{Re}}\,\O$. Second, we use the monotonicity formula to
derive three integral estimates, which roughly say that when
averaged over any time interval the mean curvature vector
$\h2_\lambda$ and the phase function $\cos\theta_\lambda$ both
tend to 0 in $L^2$ norm over a fixed ball near the singularity, as
$\lambda\to \infty$. Another direct consequence of the
monotonicity formula is that there is an upper bound of the volume
density of the rescaled submanifolds $\S^\lambda_t$, which allows
us to extract converging subsequence in measure.

\subsection{A weighted monotonicity formula}

Let $H (\x2,\x2_0,t_0,t)$ be the backward heat kernel on ${\bf
R}^k$. Let $N_t$ be a smooth family of submanifolds of dimension
$n$ in ${\bf R}^k$ defined by $F_t:N\rightarrow {\bf R}^k$. Define
\begin{eqnarray}\label{general rho}
\rho (\x2,t)&=&(4\pi (t_0-t))^{(k-n)/2}H (\x2,\x2_0,t_0,t)\nonumber\\
&=&\fr{1}{(4\pi (t_0-t))^{n/2}}\exp \left(-\fr{|\x2-\x2_0|^2}
{4(t_0-t)}\right)
\end{eqnarray}
for $t<t_0$.

A straightforward calculation (cf. [CL1], [H1], [Wa]) shows
$$
\frac{\partial}{\partial t}\rho= \left(
\frac{n}{2(t_0-t)}-\frac{\h2\cdot (\x2-\x2_0)}{2(t_0-t)}
-\frac{|\x2-\x2_0|^2}{4(t_0-t)^2}\right)\rho
$$
and along $N_t$
$$
\btu\rho=
\left(\fr{|\la\x2-\x2_0,\btd\x2\ra |^2}{4(t_0-t)^2}\right.
\left. -\fr{\la \x2-\x2_0,\btu\x2\ra}{2(t_0-t)}
-\fr{|\btd\x2|^2}{2(t_0-t)}\right)\rho
$$
where $\Delta,\nabla$ are on $N_t$ in the induced metric. Let
$N_t=\S_t$ be a smooth 1-parameter family of compact Lagrangian
submanifolds in a compact Calabi-Yau manifold $(M,\O)$ of complex
dimension $n$. Note that in the induced metric on $\Sigma_t$
$$
|\btd\f2 |^2=n ~~ \hbox{and} ~~ \btu \f2=\h2.
$$
Therefore
\begin{equation}\label{rho}
\left(\frac{\partial}{\partial t}+\Delta\right)\rho=
-\left(\left|\h2+\frac{(\f2-\x2_0)^\perp}{2(t_0-t)}\right|^2-|\h2|^2\right)\rho.
\end{equation}

On $\S_t$ we set
$$
v=\cos\theta.
$$
Denote the injectivity radius of $(M,g)$ by $i_M$. For $\x2_0\in
M$, take a normal coordinate neighborhood $U$ and let $\phi \in
C^\infty_0(B_{2r}(\x2_0))$ be a cut-off function with $\phi\equiv
1$ in $B_r(\x2_0))$, $0<2r<i_M$. Using the local coordinates in
$U$ we may regard $F(x,t)$ as a point in ${\bf R}^{2n}$ whenever
$F(x,t)$ lies in $U$. We define
$$
\Psi(\x2_0,t_0,t)=\int_{\S_t}\frac{1}{v}\phi(F)\rho(F,\x2_0,t,t_0)d\mu_t
$$
where $\rho$ is defined by (\ref{general rho}) by taking $k=2n$.

\begin{proposition}\label{mono}
Let $F_t:\S\to M$ be a smooth mean curvature flow of a compact Lagrangian submanifold
$\S_0$ in a compact Calabi-Yau manifold $M$ of complex dimension $n$. Suppose that $\S_0$
is almost calibrated by $\mbox{Re}\O$.
Then there are positive
constants $c_1$ and $c_2$ depending only on $M$, $\f2_0$ and
$r$ which is the constant in the definition of $\phi$, such that
\begin{eqnarray} \lefteqn{
\fr{\p}{\p t}\left(e^{c_1\sqrt{t_0-t}}\int_{\g1_t}\fr{1}{v}\phi\rho d\mu_t\right)}
\nonumber\\
&\leq&-e^{c_1\sqrt{t_0-t}}\int_{\S_t}\frac{1}{v}\phi\rho\left(\frac{2|\nabla
v|^2}{v^2}
+\left|\h2+\fr{(\f2-\x2_0)^{\perp}}{2(t_0-t)}\right|^2+\frac{|\h2|^2}{2}\right)\no\\
&& +c_2e^{c_1\sqrt{t_0-t}}. \label{mon2}
\end{eqnarray}
\end{proposition}
{\it Proof.} Notice that
$$
\btu F=\h2+g^{ij}\Gamma^\alpha_{ij}v_\alpha
$$
where $v_\alpha, \alpha=1,...,n$ is a basis of $T^\perp\S_t$, $g^{ij}$ is the induced metric on $\S_t$
and $\Gamma^\alpha_{ij}$ is the Christoffel symbol on $M$. Equation (\ref{rho}) reads as
\begin{equation}
\left(\frac{\partial}{\partial t}+\Delta\right)\rho=
-\left(\left|\h2+\frac{(\f2-\x2_0)^\perp}{2(t_0-t)}\right|^2-|\h2|^2
+\frac{g^{ij}\Gamma^\alpha_{ij}v_\alpha\cdot(\f2-\x2_0)}{t_0-t}\right)\rho.
\end{equation}
>From  (\ref{cos}) we have
$$
\frac{\partial}{\partial t}\frac{1}{v}=
\Delta\frac{1}{v}-\frac{|\h2|^2}{v}-\frac{2|\nabla v|^2}{v^3}
$$
and
$$
\frac{d}{dt}d\mu_t=-|\h2|^2d\mu_t.
$$
Moreover
$$
\frac{\p\phi(F)}{\p t}=\nabla\phi\cdot\h2.
$$
Now we have
\begin{eqnarray}\label{estimate3} \lefteqn{
\frac{d}{dt}\int_{\S_t}\frac{1}{v}\phi\rho }\nonumber\\
&=&\int_{\S_t}\phi\rho\Delta \frac{1}{v}
-\int_{\S_t}\left(\frac{|\h2|^2}{v}+\frac{2}{v^3}|\nabla v|^2\right)\phi\rho
+\int_{\S_t}\frac{1}{v}\nabla \phi \cdot\h2 \rho \nonumber\\
&&-\int_{\S_t}\frac{1}{v}\phi\left(\Delta\rho
+\left(\left|\h2+\frac{(\f2-\x2_0)^\perp}{2(t_0-t)}\right|^2-|\h2|^2
+\frac{g^{ij}\Gamma^\alpha_{ij}v_\alpha\cdot(\f2-\x2_0)}{t_0-t}\right)\rho\right)\nonumber\\
&&-\int_{\S_t}\frac{1}{v}\phi\rho|\h2|^2\nonumber\\
&\leq&-\int_{\S_t}\phi\rho\left(\frac{2}{v^3}|\nabla v|^2
+\frac{1}{v}\left|\h2+\frac{(\f2-\x2_0)^\perp}{2(t_0-t)}\right|^2
+\frac{|\h2|^2}{v}\right)\nonumber\\
&&+\int_{\S_t}\left(\phi\rho\Delta\frac{1}{v}-\frac{1}{v}\phi\Delta\rho\right)
-\int_{\S_t}\frac{1}{v}\phi\rho\frac{g^{ij}
\Gamma^\alpha_{ij}v_\alpha\cdot(\f2-\x2_0)}{t_0-t}\nonumber\\
&&+\int_{\S_t}\frac{1}{v}\rho\left(\epsilon^2 \phi |\h2|^2
+\frac{1}{4\epsilon^2}\frac{|\nabla\phi|^2}{\phi}\right)
\end{eqnarray}
where we used Cauchy-Schwartz inequality for $\nabla\phi\cdot\h2$. By Stokes formula
$$
\int_{\S_t}\left(\phi\rho\Delta\frac{1}{v}-\frac{1}{v}\phi\Delta\rho\right)
=2\int_{\S_t}\frac{1}{v}\nabla\phi\nabla\rho+\int_{\S_t}\frac{1}{v}\rho\Delta\phi.
$$
Since $\phi\in C^\infty_0(B_{2r}(\x2_0),{\bf R}^+)$, we have (cf. Lemma 6.6 in [Il])
$$
\frac{|\nabla\phi|^2}{\phi}\leq 2\max_{\phi>0}|\nabla^2\phi|.
$$
Note that $\nabla\phi\equiv 0$ in $B_r(\x2_0)$, so $|\rho\Delta\phi|$ and 
$|\nabla\phi\cdot\nabla\rho|$ are bounded in $B_{2r}(\x2_0)$. Hence
\begin{equation}\label{estimate2}
\int_{\g1_t}\left|\frac{1}{v}\rho \btu \phi\right| +
\int_{\g1_t}\left|\frac{1}{v}\btd \phi\cdot\btd\rho \right|\leq
C\int_{\g1_t}\frac{1}{v}d\mu_t\leq \frac{C}{\min_{\S_0}v}\mbox{vol}(\S_0)
\end{equation}
where $C$ depends only on $r,\max(|\nabla^2\phi|+|\nabla\phi|)$.

Since $\Gamma^\al_{ij} (\x2_0)=0$, we may choose $r$ sufficiently
small such that
$$
|g^{ij}\Gamma^\al_{ij} (\f2 )|\leq C|\f2 -\x2_0|
$$
in $B_{2r}(\x2_0)$ for some constant $C$ depending on $M$.
We claim
\begin{equation}\label{estimate1}
\fr{|g^{ij}\Gamma^\al_{ij} v_{\al}\cdot(F-\x2_0)|}{t_0-t}\rho
(\f2 ,t)\leq c_1\fr{\rho (\f2 ,t)}{\sqrt{t_0-t}}+C.
\end{equation}
In fact it suffices to show for any $x$ and $s>0$
$$
\frac{x^2}{s}\frac{e^{-x^2/s}}{s^{n/2}}\leq
C\left(1+\frac{1}{s^{1/2}}\frac{e^{-x^2/s}}{s^{n/2}}\right).
$$
To see this, let $y=x^2/s$ and then it is easy to verify that
$$
y\leq C\left(s^{n/2}e^y+\frac{1}{s^{1/2}}\right)
$$
holds trivially if $y\leq1/s^{1/2}$ and follows from $y^{n+1}\leq Ce^y$ if
$y>1/s^{1/2}$ for some $C$. So (\ref{estimate1}) is established.

Letting $\epsilon^2=1/2$ in (\ref{estimate3}) and applying
(\ref{estimate2}), (\ref{estimate1}) to (\ref{estimate3})
we have
$$
\fr{\p}{\p t}\Psi\leq -\int_{\g1_t}\frac{1}{v}\phi\rho\left(\frac{2|\nabla v|^2}{v^2}
+\left|\h2+\fr{(\f2-\x2_0)^{\perp}}{2(t_0-t)}\right|^2+\frac{|\h2|^2}{2}\right)
+\fr{c_1}{\sqrt{t_0-t}}\Psi +c_2.
$$
The proposition follows. \hfill Q.E.D.

Suppose that $(X_0,T)$ is a singular point of the mean curvature
flow (\ref{meaneqn}). We now describe the rescaling process around
$(X_0,T)$. For any $t<0$, we set
\begin{equation}\label{def:blowup}
F_\lmd(x,t)=\lmd (F(x,T+\lmd^{-2}t)-X_0)
\end{equation}
where $\lmd$ are positive constants which go to infinity. The
scaled submanifold is denoted by $\g1_t^\lmd=F_\lmd(\Sigma,t)$ on
which $d\mu^\lmd_t$ is the area element obtained from $d\mu_t$. If
$g^\lmd$ is the metric on $\Sigma^\lmd_t$, it is clear that
$$
g_{ij}^\lmd=\lmd^2g_{ij},~~~~(g^\lmd)^{ij}=\lmd^{-2}g^{ij}.
$$
We therefore have
\begin{eqnarray*}
\fr{\p\f2_\lmd}{\p t}&=&\lmd^{-1}\fr{\p\f2}{\p t}\\
\h2_\lmd &=&\lmd^{-1}\h2\\
|\a2_\lmd|^2&=&\lmd^{-2}|\a2|^2.
\end{eqnarray*}
It follows that the scaled submanifold also evolves by a mean
curvature flow
\begin{equation}
\fr{\p\f2_\lmd}{\p t}=\h2_\lmd .
\end{equation}
Moreover, since
\begin{eqnarray*}
d\mu^\lmd_t(F_\lmd(x,t))&=&\lmd^n d\mu_t(F(x,T+\lmd^{-2}t))\\
\O|_{\S^\lmd_t}(F_\lmd(x,t))&=&\lmd^n\O|_{\S_t}(F(x,T+\lmd^{-2}t))
\end{eqnarray*}
we have
$$
\cos\theta_\lmd (F_\lmd(x,t))=\cos\theta (F(x,T+\lmd^{-2}t)).
$$

\subsection{Integral estimates}

\begin{proposition}\label{main2.1}
Let $(M,\O)$ be a Calabi-Yau manifold of complex dimension $n$. If the
initial compact submanifold is Lagrangian and is almost calibrated by ${\mbox Re}\O$,
then for any $R>0$ and any
$-\infty<s_1<s_2<0$, we have
\be\label{m2.2}
\int_{s_1}^{s_2}\int_{\g1_t^\lmd\cap B_R(0)} |\btd \cos
\theta_{\lmd}|^2d\mu_t^\lmd dt\to 0~~{\rm as}~~\lmd\to\infty, \ee
\be\label{m2.3} \int_{s_1}^{s_2}\int_{\g1_t^\lmd\cap
B_R(0)}|\h2_\lmd|^2d\mu_t^\lmd dt\to 0~~{\rm as}~~\lmd\to\infty,
\ee and \be\label{m2.4} \int_{s_1}^{s_2}\int_{\g1_t^\lmd\cap
B_R(0)}|F_\lmd^{\perp}|^2d\mu_t^\lmd dt\to 0~~{\rm
as}~~\lmd\to\infty. \ee
\end{proposition}
{\it Proof:} For any $R>0$, we choose a cut-off function
$\phi_R\in C_0^{\infty}(B_{2R}(0))$ with $\phi_R\equiv 1$ in
$B_R(0)$, where $B_r(0)$ is the metric ball centered at $0$ with
radius $r$ in ${\bf R}^{2n}$. For any fixed $t<0$, the mean curvature flow (\ref{meaneqn}) has a
smooth solution near $T+\lmd^{-2}t<T$ for sufficiently large $\lmd$, since
$T>0$ is the first blow-up time of the flow.
Let $v_\lmd=\cos\theta_\lmd$. It is clear
\bs
\lefteqn{\int_{\g1_t^\lmd}\fr{1}{v_\lmd}\fr{1}{(0-t)^{n/2}}\phi_R(\f2_\lmd
)\exp
\left(-\fr{|\f2_\lmd |^2}{4(0-t)}\right)d\mu_t^\lmd}\\
&=&\int_{\g1_{T+\lmd^{-2}t}}\fr{1}{v_\lmd}\phi (\f2_\lmd)
\fr{1}{(T-(T+\lmd^{-2}t))^{n/2}}\exp
\left(-\fr{|\f2(x,T+\lmd^{-2}t)-X_0|^2}{4(T-(T+\lmd^{-2}t))}\right)d\mu_t
, \es where $\phi$ is the function defined in the definition of
$\Phi$. Note that $T+\lmd^{-2}t\to T$ for any fixed $t$ as $\lmd\to\infty$. By
the weighted monotonicity formula (\ref{mon2}),
$$
\frac{\partial}{\partial t}\left(e^{c_1\sqrt{t_0-t}}\Psi\right)
\leq c_2e^{c_1\sqrt{t_0-t}},
$$
and it then follows that $\lim_{t\to t_0}e^{c_1\sqrt{t_0-t}}\Psi$
exists. This implies, by taking $t_0=T$ and $t=T+\lmd^{-2}s$, that
for any fixed $s_1$ and $s_2$ with $-\infty <s_1<s_2<0$,
\begin{eqnarray}\label{limit}
\lefteqn{e^{c_1\sqrt{T-(T+\lmd^{-2}s_2)}}
\int_{\g1_{s_2}^\lmd}\fr{1}{v_\lmd}\phi_R\fr{1}{(0-s_2)^{n/2}}\exp
\left(-\fr{|\f2_\lmd|^2}{4(0-s_2)}\right)d\mu_{s_2}^\lmd}\nonumber\\
&&-e^{c_1\sqrt{T-(T+\lmd^{-2}s_1)}}\int_{\g1_{s_1}^\lmd}\fr{1}{v_\lmd}
\phi_R\fr{1}{(0-s_1)^{n/2}}\exp
\left(-\fr{|\f2_\lmd |^2}{4(0-s_1)}\right)d\mu_{s_1}^\lmd \nonumber\\
& \to &0 ~~\hbox{as}~~ \lmd\to\infty.
\end{eqnarray}
Integrating (\ref{mon2}) from $T+\lmd^{-2}s_1$ to
$T+\lmd^{-2}s_2$, and letting $T+\lmd^{-2}s=t$, we get
\ba\label{djh2}
\lefteqn{-e^{c_1\sqrt{-\lmd^{-2}s_2}}\int_{\g1_{s_2}^\lmd}\fr{1}{v_\lmd}
\phi_R\fr{1}{(0-s_2)^{n/2}}\exp
\left(-\fr{|\f2_\lmd|^2}{4(0-s_2)}\right)d\mu_{s_2}^\lmd}\no\\
&&+e^{c_1\sqrt{-\lmd^{-2}s_1}}\int_{\g1_{s_1}^\lmd}\fr{1}{v_\lmd}\phi_R\fr{1}{(0-s_1)^{n/2}}\exp
\left(-\fr{|\f2_\lmd|^2}{4(0-s_1)}\right)d\mu_{s_1}^\lmd\no \\
&\geq
&\int_{T+\lmd^{-2}s_1}^{T+\lmd^{-2}s_2}e^{c_1\sqrt{T-t}}\int_{\S_t}\frac{1}{v}\phi\rho\left(\frac{2|\nabla
v|^2}{v^2}
+\left|\h2+\fr{(\f2-\x2_0)^{\perp}}{2(T-t)}\right|^2+\frac{|\h2|^2}{2}\right)d\mu_t\no\\
&&-c_2\lmd^{-2}(s_2-s_1)\no\\
&\geq &
\int_{s_1}^{s_2}e^{c_1\sqrt{-\lmd^{-2}s}}\int_{\g1_s^\lmd}\fr{1}{v_\lmd}\phi_R\rho
(\f2_\lmd ,s)\left|\h2_\lmd+\fr{(\f2_\lmd)^{\perp}}{2(-s)}\right|^2d\mu_s^\lmd\no\\
&&+\int_{s_1}^{s_2}e^{c_1\sqrt{-\lmd^{-2}s}}\int_{\g1_s^\lmd}\fr{1}{v_\lmd}\phi_R\rho
(\f2_\lmd,s)\frac{|\h2_\lmd|^2}{2}d\mu_s^\lmd\no\\
&&+\int_{s_1}^{s_2}e^{c_1\sqrt{-\lmd^{-2}s}}\int_{\g1_s^\lmd}\fr{2}{v_\lmd^3}|\btd
v_\lmd|^2\phi_R\rho (\f2_\lmd
,s)d\mu_s^\lmd\no\\
&&-c_2\lmd^{-2}(s_2-s_1). \ea From (\ref{limit}) and (\ref{djh2})
the proposition follows.   \hfill Q.E.D.

\subsection{Upper bound on volume density}
Now we show the existence of the $\lmd$ tangent cones by deriving an finite upper bound for
the volume density. These cones are independent of $t$, but may depend on the blowing up sequence
$\lmd$.

\begin{proposition}\label{exit}
Suppose that $\Sigma_t$ evolves along mean curvature flow and $\S_0$ is a compact Lagrangian 
submanifold in $(M,\O)$ and is almost calibrated by $\mbox{Re}\,\O$. 
For any $\lmd,R>0$ and any $t<0$, \be\label{fm}
\mu_t^{\lmd}(\g1_t^\lmd\cap B_R(0))\leq CR^n,
 \ee
where $B_R(0)$ is a metric ball in ${\bf R}^{2n}$ and $C>0$ is
independent of $\lmd$. For any sequence $\lmd_i\to\infty$, there
is a subsequence $\lmd_k\to\infty$ such that
$(\g1_t^{\lmd_k},\mu_t^{\lmd_k})\to (\g1^\infty,\mu^\infty)$ in
the sense of measure, for any fixed $t<0$, where
$(\g1^\infty,\mu^\infty)$ is independent of $t$. The multiplicity
of $\g1^\infty $ is finite.
\end{proposition}
{\it Proof:}  We shall first prove the inequality (\ref{fm}). We shall
use $C$ below for uniform positive constants which are independent
of $R$ and $\lmd$. Straightforward computation shows
\bs
\mu_t^{\lmd}(\g1_t^\lmd\cap B_R(0))&=
&\lmd^n\int_{\g1_{T+\lmd^{-2}t}\cap B_{\lmd^{-1}R}(X_0)}d\mu_t\\
&=& R^n(\lmd^{-1}R)^{-n}\int_{\g1_{T+\lmd^{-2}t}\cap B_{\lmd^{-1}R}(X_0)}d\mu_t\\
&\leq & CR^n\int_{\g1_{T+\lmd^{-2}t}\cap B_{\lmd^{-1}R}(X_0)}\frac{1}{v_\lmd}
\frac{1}{(4\pi)^{n/2}(\lmd^{-1}R)^n}e^{-\frac{|X-X_0|^2}{4(\lmd^{-1}R)^2}}d\mu_t\\
&=&CR^n\Psi(X_0,T+(\lmd^{-1}R)^2+\lmd^{-2}t, T+\lmd^{-2}t).
\es
By the weighted monotonicity inequality (\ref{mon2}), we have
\bs
\mu_t^{\lmd}(\g1_t^\lmd\cap B_R(0))&\leq &CR^n\Psi(X_0,
T+(\lmd^{-1}R)^2+\lmd^{-2}t, T/2)+CR^n\\
&\leq &\fr{\mu_{T/2}(\g1_{T/2})}{T^{n/2}\min_{\S_0}v}CR^n+CR^n.
\es
Since volume is non-increasing along mean curvature flow:
$$
\fr{\p}{\p t}\mu_t(\g1_t)=-\int_{\g1_t}|\h2|^2d\mu_t,
$$
we have therefore established (\ref{fm}):
$$
\mu_t^{\lmd}(\g1_t^\lmd\cap B_R(0))\leq CR^n.
$$
By (\ref{fm}), the compactness theorem of the measures (c.f.
[Si1], 4.4) and a diagonal subsequence argument, we conclude
that there is a subsequence $\lmd_k\to\infty$ such that
$(\g1_{t_0}^{\lmd_k},\mu_{t_0}^{\lmd_k})\to
(\g1_{t_0}^\infty,\mu_{t_0}^\infty)$ in the sense of measures for a
fixed $t_0<0$.

We now show that, for any $t<0$, the subsequence $\lmd_k$ which we
have chosen above satisfies
$(\g1_{t}^{\lmd_k},\mu_{t}^{\lmd_k})\to
(\g1_{t_0}^\infty,\mu_{t_0}^\infty)$ in the sense of measure. And
consequently the limiting submanifold
$(\g1^\infty_{t_0},\mu^\infty_{t_0})$ is independent of $t_0$.
Recall that the following standard formula for mean curvature
flows
\begin{equation}\label{variation}
\frac{d}{dt}\int_{\Sigma^\lmd_t}\phi d\mu^\lmd_t
=-\int_{\Sigma^\lmd_t}\left(\phi|\h2_\lmd
|^2+\nabla\phi\cdot\h2_\lmd\right)d\mu^\lmd_t
\end{equation}
is valid for any test function $\phi\in C_0^\infty(M)$ (cf. (1) in
Section 6 in [I2] and [B] in the varifold setting).

Then for any given $t<0$ integrating (\ref{variation}) yields
 \bs
 \int_{\Sigma^{\lmd_k}_t}\phi d\mu^{\lmd_k}_t
 -\int_{\Sigma^{\lmd_k}_{t_0}}\phi d\mu^{\lmd_k}_{t_0} &=&
 \int_t^{t_0}\int_{\Sigma^{\lmd_k}_t}\left(\phi|\h2_{\lmd_k}
|^2+\nabla\phi\cdot\h2_{\lmd_k}\right)d\mu^{\lmd_k}_tdt\\
&\to & 0~~{\rm as}~~k\to\infty~~{\rm by}~~(\ref{m2.3}) .\es So,
for any fixed $t<0$, $(\g1_{t}^{\lmd_k},\mu_{t}^{\lmd_k})\to
(\g1_{t_0}^\infty,\mu_{t_0}^\infty)$ in the sense of measures as
$k\to\infty$. We denote $(\g1_{t_0}^\infty,\mu_{t_0}^\infty)$ by
$(\g1^\infty,\mu^\infty)$, which is independent of $t_0$.

The inequality (\ref{fm}) yields a uniform upper bound on
$R^{-n}\mu^{\lmd_k}_t(\g1^{\lmd_k}_t\cap B_R(0))$, which
yields finiteness of the multiplicity of $\g1^\infty$.
\hfill Q.E.D.

\begin{definition}\label{bubble1}\label{bb}
{\em Let $(X_0,T)$ be a singular point of the mean curvature flow
of a compact Lagrangian submanifold $\Sigma_0$ in a compact
Calabi-Yau manifold $M$. We call $(\g1^\infty, d\mu^\infty)$
obtained in Proposition \ref{exit} {\it a $\lmd$ tangent cone of
the mean curvature flow $\g1_t$ at $(X_0,T)$.}}
\end{definition}

\section{Rectifiability of $\lmd$ tangent cones}

In this section we shall show that the $\lmd$ tangent cone $\g1^\infty$ is ${\mathcal
H}^n$-rectifiable, where ${\mathcal H}^n$ is the $n$-dimensional Hausdorff measure.

\begin{proposition}\label{rectif}
Let $M$ be a compact Calabi-Yau manifold of complex dimension $n$. If the initial compact
submanifold $\g1_0$ is Lagrangian and almost calibrated by ${\mbox{Re}}\,\O$,
then the $\lmd$ tangent cone $(\g1^\infty
,d\mu^\infty )$ of the mean curvature flow at $(X_0,T)$ is
${\mathcal H}^n$-rectifiable.
\end{proposition}
{\it Proof:} Let
$(\S_t^k,d\mu_t^k)=(\S_t^{\lmd_k},d\mu_t^{\lmd_k})$. We set
$$
A_R=\left\{t\in (-\infty,0)\left| ~\lim_{k\to
\infty}\int_{\g1_t^k\cap B_R(0)}|{\bf H}_k|^2d\mu_t^k
\not=0\right\}, \right.
$$
and
$$
A=\bigcup_{R>0}A_R.
$$

Denote the measures of $A_R$ and $A$ by $|A_R|$ and $|A|$
respectively. It is clear from (\ref{m2.3}) that $|A_R|=0$ for any $R>0$.
So $|A|=0$.

For any $\xi\in\g1^\infty$, choose $\xi_k\in\g1_t^k$
with $\xi_k\to\xi$ as $k\to\infty$. By the monotonicity identity
(17.4) in [Si1], we have
\ba\label{smon}
\sigma^{-n}\mu_t^k(B_\sigma(\xi_k))&=&\rho^{-n}\mu_t^k(B_\rho(\xi_k))
-\int_{B_\rho(\xi_k)\setminus B_\sigma(\xi_k)}\fr{|D^\perp r|^2}{r^n}d\mu_t^k\no\\
&& -\fr{1}{n}\int_{B_\rho(\xi_k)}(x-\xi_k)\cdot {\bf H}_k
\left(\fr{1}{r_\sigma^n}-\fr{1}{\rho^n}\right)d\mu_t^k, \ea for
all $0<\sigma \leq \rho$, where $\mu_t^k(B_\sigma(\xi_k))$ is the
measure of $\Sigma^k_t\cap B_\sigma((\xi_k))$, $r=r(x)$ is the
distance from $\xi_k$ to $x$, $r_\sigma =\max \{r,\sigma\}$, and
$D^\perp r$ denotes the orthogonal projection of $Dr$ (which is a
vector of length 1) onto $(T_{\xi_k}\g1^k_t)^\perp$. Choosing
$t\not\in A$, we have
$$
\lim_{k\to\infty}\int_{B_\rho(\xi_k)}|{\bf H}_k|^2d\mu_t^k=0.
$$
H\"older's inequality and (\ref{fm}) then lead to
\begin{eqnarray}\label{estimate}
\lefteqn{\lim_{k\to\infty}\left|\int_{B_\rho(\xi_k)}(x-\xi_k)\cdot {\bf H}_k
\left(\fr{1}{r_\sigma^n}-\fr{1}{\rho^n}\right)d\mu_t^k\right|}\nonumber\\
&\leq&
C\rho\left(\fr{1}{\sigma^n}-\fr{1}{\rho^n}\right)\lim_{k\to\infty}
\left(\sqrt{\mu^k_t(B_{\rho}(\xi_k))}
\sqrt{\int_{B_\rho(\xi_k)}|\h2_k|^2 d\mu^k_t}\right)\nonumber\\
&\leq& C\rho^{1+n/2}\left(\fr{1}{\sigma^n}-\fr{1}{\rho^n}\right)
\lim_{k\to\infty}\sqrt{\int_{B_\rho(\xi_k)}|\h2_k|^2 d\mu^k_t}\nonumber\\
&=&0.
\end{eqnarray}
Letting $k\to\infty$ in (\ref{smon}) and using (\ref{estimate}), we obtain
$$
\sigma^{-n}\mu^\infty(B_\sigma(\xi))\leq
\rho^{-n}\mu^\infty(B_\rho(\xi)),
$$
for all $0<\sigma \leq \rho$. By (\ref{fm}) we know that
$$
 \lim_{\rho\to 0}\rho^{-n}\mu^\infty(B_\rho(\xi))<C<\infty.
$$
Therefore, $\lim_{\rho\to 0}\rho^{-n}\mu^\infty(B_\rho(\xi))$
exists.

We shall show that there exists a positive number $r_0$ such that
for any $0<r<r_0<1$ the following density estimate holds
\be\label{rec}
\lim_{\rho\to 0}\rho^{-n}\mu^\infty(B_\rho(\xi))\geq \fr{1}{4c(n)+4}>0 \ee
for some positive constant $c(n)$ which will be determined below.
Assume (\ref{rec}) fails to hold. Then there is $\rho_0>0$ such
that
$$
(2\rho_0)^{-n}\mu^\infty(B_{2\rho_0}(\xi))<\frac{1}{4c(n)+4}.
$$
By the monotonicity formula (\ref{smon}) and that $\mu^k_t$
converges to $\mu^\infty$ as measures, there exists $k_0>0$ such
that, for all $0<\rho<2\rho_0$ and $k>k_0$, we have
\be\label{rec1} \rho^{-n}\mu_t^k(B_\rho(\xi))< \fr{1}{2c(n)+2}.\ee
Take a cut-off function $\phi_\rho\in C^{\infty}_0(B_\rho(\xi))$
on the $2n$-dimensional ball $B_\rho(\xi_k)$ so that
\begin{eqnarray*}
&&\phi_\rho\equiv 1~~{\rm in}~~B_{\fr{\rho}{2}}(\xi)\\
&&0\leq\phi_\rho\leq 1,~~{\rm and}~~| \btd\phi_\rho
|\leq\fr{C}{\rho},~~{\rm in}~ B_\rho(\xi).
\end{eqnarray*}
>From (\ref{variation}), we have \bs
\lefteqn{\rho^{-n}\int_{B_\rho(\xi)}\phi_\rho d\mu_{t-r^2}^k
-\rho^{-n}\int_{B_\rho(\xi)}\phi_\rho d\mu_{t}^k
}\\
&\leq & C\rho^{-n}\int_{t-r^2}^t\int_{B_{\rho}(\xi)}|{\bf H}_k|
^2d\mu_s^kds+C\rho^{-n-1}\int_{t-r^2}^t\int_{B_\rho(\xi)}|{\bf H}_k|
d\mu_s^kds\\
&\leq & C\rho^{-n}\int_{t-r^2}^t\int_{B_{\rho}(\xi)}|{\bf H}_k|
^2d\mu_s^kds +
C\rho^{-n-1}\int_{t-r^2}^t\left(\int_{B_{\rho}(\xi)}|{\bf H}_k|
^2d\mu_s^k\right)^{1/2}\mu_s^k(B_\rho(\xi))^{1/2}ds \\
&\leq &
C\rho^{-n}\int_{t-r^2}^t\int_{B_{\rho}(\xi)}|{\bf H}_k|
^2d\mu_s^kds +
C\rho^{-n/2-1}\int_{t-r^2}^t\left(\int_{B_{\rho}(\xi)}|{\bf H}_k|
^2d\mu_s^k\right)^{1/2}ds~~{\rm by}~(\ref{fm})\\
&\to& 0, ~~{\rm as}~~ k\to\infty~~{\rm by}~(\ref{m2.3}). \es Here
we have used $C$ for uniform positive constants which are
independent of $k$ and $\rho$.
Therefore, there are constants $\d1_1>0$ and $k_1>0$ such that for
all $\rho$ and $k$ with $0<\rho<\d1_1$, $0<r<1$, and $k>k_1$ the
estimate
\be \label{rec2}
\rho^{-n}\mu_{t-r^2}^k(B_\rho(\xi))<\frac{1}{c(n)+1}< 1 \ee holds. Let
$d\sigma^k_{t-r^2}$ be the area element of $\partial B_{\rho}(\xi)\cap \g1^k_{t-r^2}$.
By the co-area formula, for $0<r<<1$, for a smooth cut-off function $\phi$ with
support in the $2n$-dimensional ball $B_{\delta_1}(0)$ in ${\bf R}^{2n}$ with
$0\leq\phi\leq 1, \phi\equiv 1$ in $B_{\delta_1/2}(0)$, we have
\begin{eqnarray}\label{density}
\Phi_k(\xi,t,t-r^2) &=&\frac{1}{(4\pi r^2)^{n/2}}\int_{\Sigma_{t-r^2}^k}\phi\,
e^{-\frac{|F_k-\xi|^2}{4r^2}}d\mu^k_{t-r^2}\nonumber\\
&\leq &\fr{1}{(4\pi)^{n/2} r^n}\int_0^{\d1_1}\int_{\part
 B_\rho(\xi)\cap\g1_{t-r^2}^k}e^{-\fr{\rho^2}{4r^2}}d\sigma_{t-r^2}^k d\rho\nonumber\\
&\leq &\fr{1}{(4\pi)^{n/2} r^n}\int_0^{\d1_1}
e^{-\fr{\rho^2}{4r^2}}\frac{d}{d\rho}{\rm
Vol}(B_\rho(\xi)\cap\Sigma^k_{t-r^2})d\rho\nonumber\\
&\leq&\frac{1}{\pi^{n/2}(2r)^n}{\mbox{Vol}}(B_{\delta_1}(\xi)\cap\S^k_{t-r^2}))
e^{-\frac{\delta^2_1}{4r^2}}\nonumber\\
&&+\frac{1}{(c(n)+1)\pi^{n/2}}\int^{\delta_1}_0
e^{-\frac{\rho^2}{4r^2}}\frac{\rho^n}{2^nr^n}d\frac{\rho^2}{4r^2}\nonumber
\end{eqnarray}
by integration by parts and (\ref{rec2}). By (\ref{fm}),
$$
\frac{1}{\pi^{n/2}(2r)^n}{\mbox{Vol}}(B_{\delta_1}(\xi)\cap\S^k_{t-r^2}))
e^{-\frac{\delta^2_1}{4r^2}}\leq C\left(\frac{\delta_1}{2r}\right)^n
e^{-\frac{\delta^2_1}{4r^2}}=o(r).
$$
Letting $y=\rho/2r$ we have
\begin{eqnarray*}
\int_0^{\delta_1}e^{-\frac{\rho^2}{4r^2}}\left(\frac{\rho}{2r}\right)^n
d\left(\frac{\rho}{2r}\right)^2\leq 2\int^\infty_0 e^{-y^2}y^{n+1}dy=c(n)<\infty,
\end{eqnarray*}
and there is an explicit formula for $c(n)$ depends on whether $n$ is odd or even.
Thus we conclude
$$
\Phi_k(\xi,t,t-r^2)\leq 1+o(r).
$$
For any classical mean curvature flow $\Gamma_t$ in a compact
Riemannian manifold which is isometrically embedded in ${\bf
R}^N$, White proves a local regularity theorem (Theorem 3.1 and
Theorem 4.1 in [Wh1]): When dim$\Gamma_t=n$,
 there is a constant $\epsilon>0$ such that if the
Gaussian density satisfies
$$
\lim_{r\rightarrow 0}\int_{\Gamma_{t-r^2}}\frac{1}{(4\pi r^2)^{n/2}}
\exp\left(-\fr{|y-x|^2}{4r^2}\right)d\mu(y) <1+\epsilon,
$$
then the mean curvature flow is smooth in a neighborhood of $x$.
Combining this regularity result with (\ref{density}), we are led
to choose $r>0$ sufficiently small and then conclude that
$$
\sup_{B_r(\xi)\cap\g1_t^k}|\a2_k|\leq C
$$
and consequently $\g1_t^k$ converges strongly in $B_r(\xi)\cap\g1_t^k$ to
$\Sigma^\infty_t\cap B_r(\xi)$, as $k\to\infty$. So $\Sigma^\infty\cap B_r(\xi)$ is
smooth. Smoothness of $\Sigma^\infty\cap B_r(\xi)$ immediately implies
$$
\lim_{\rho\to 0}\rho^{-n}\mu^\infty(B_\rho(\xi))= 1.
$$
This contradicts (\ref{rec1}). Hence we have established
(\ref{rec}).

In summary, we have shown that $\lim_{\rho\to
0}\rho^{-n}\mu^\infty(B_\rho(\xi))$ exists and for ${\mathcal
H}^n$ almost all $\xi\in\g1^\infty$,
\begin{equation}\label{den}
\fr{1}{4c(n)+4}\leq\lim_{\rho\to
0}\rho^{-n}\mu^\infty(B_\rho(\xi))<\infty.
\end{equation}

Finally, we recall a fundamental theorem of Priess in [P]: if
$0\leq m\leq p$ are integers and $\Omega$ is a Borel measure on
${\bf R}^p$ such that
$$
0<\lim_{r\rightarrow 0}\fr{\Omega(B_r(x))}{r^m}<\infty,
$$
for almost all $x\in\Omega$, then $\Omega$ is $m$-rectifiable. Now
we conclude from (\ref{den}) that $(\g1^\infty, \mu^\infty)$ is
${\mathcal H}^n$-rectifiable. \hfill
Q.E.D.

\section{Minimality of the $\lmd$ tangent cones}

In this section, we will show that the $\lmd$ tangent cone
$\g1^\infty$ is a stationary integer multiplicity rectifiable current in ${\bf R}^{2n}$.

\begin{theorem}\label{stat}
Let $M$ be a compact Calabi-Yau manifold. If the initial compact
submanifold is Lagrangian and is almost calibrated by
$\mbox{Re}\,\O$, then the $\lmd$ tangent cone $\g1^\infty$ is a
stationary rectifiable Lagrangian current in ${\bf R}^{2n}$ with
volume density greater than one at $X_0$.
\end{theorem}
{\it Proof:} Let $V_t^k$ be the varifold defined by $\g1_t^k$. By
the definition of varifolds, we have
$$
V_t^k(\psi)=\int_{\g1_t^k}\psi(x, T{\g1_t^k})d\mu_t^k
$$
for any $\psi\in C_0^0(G^2({\bf R}^{2n}),R)$, where $G^2({\bf
R}^{2n})$ is the Grassmanian bundle of all $n$-dimensional planes
tangent to $\g1_t^\infty$ in ${\bf R}^{2n}$. For each smooth
submanifold $\g1^k_t$, the first variation $\d1V_t^k$ of $V_t^k$
(cf. [A], (39.4) in [Si1] and (1.7) in [I2]) is
$$
\d1V_t^k=-\mu_t^k\lfloor \h2_k.
$$
By Proposition \ref{main2.1}, we have that $\d1V_t^k\to 0$ at $t$
as $k\rightarrow\infty$.

Recall that a $k$-varifold is a Radon measure on $G^k(M)$,
where $G^k(M)$ is the Grassmann bundle of all $k$-planes tangent to $M$.
Allard's compactness theorem for rectifiable varifolds (6.4 in [A], also see
1.9 in [I2] and Theorem 42.7 in [Si1]) asserts the following:
let $(V_i,\mu_i)$ be a sequence of rectifiable
$k$-varifolds in $M$ with
$$
\sup_{i\geq 1}(\mu_i(U)+|\d1V_i|(U))<\infty~~{\rm
for~each}~U\subset\subset M.
$$
Then there is a varifold $(V,\mu)$ of locally bounded first variation
and a subsequence, which we
also denote by $(V_i,\mu_i)$, such that
(i)  Convergence of measures: $\mu_i\to\mu$ as Radon
measures on $M$,
(ii) Convergence of tangent planes: $V_i\to V$ as
Radon measures on $G^k(M)$,
(iii) Convergence of first variations: $\d1 V_i\to \d1 V$ as $TM$-valued
Radon measures,
(iv) Lower semicontinuity of total first variations: $|\d1
V|\leq\liminf_{i\to\infty} |\d1 V_i|$ as Radon measures.

By (iii) in Allard's compactness theorem, we have
$$
-\mu^\infty\lfloor \h2_\infty= \d1 V^\infty =\lim_{k\to\infty} \d1
V_t^k=0.
$$
Therefore $\g1^\infty$ is stationary. The rescaling process in a
neighborhood of $X_0$ in $M$ implies that the metrics $g^\lmd$
tends to the flat metric on ${\bf R}^{2n}$ and the K\"ahler 2-form
$\o^\lmd$ tends to a constant closed 2-form $\o_0$ which is
determined by $\o_0(0)=\o(X_0)$. The tangent spaces to $\S_t^k$
converge to that to $\S^\infty$ as measures by (ii) in Allard's
compactness theorem. Hence $\o^{\lmd_k}|_{\S_t^k} \to
\o_0|_{\S^\infty}$. But $\S_t^k$ is Lagrangian, it follows
$\o^{\lmd_k}|_{\S_t^k}=0$ therefore $\o_0|_{\S^\infty}=0$.
Therefore, $\S^\infty$ is a Lagrangian.

On the other hand, as $\lmd\to\infty$ in the blow-up precess, the
holomorphic $(n,0)$-form $\O$ converges to a constant holomorphic
$(n,0)$-form $\O_0$ on ${\bf R}^n$ determined by
$\O_0(0)=\O(X_0)$. We write
$$
{\mbox{Re}}\,\O_0|_{\S^\infty}=\theta_0 d\mu^\infty,
$$
$$
{\mbox{Re}}\,\O^{\lmd_k}|_{\S_t^k}=\cos\theta_{\lmd_k} d\mu_t^k
$$
and from Allard's compactness theorem 
$$
{\mbox{Re}}\,\O^{\lmd_k}|_{\S_t^k} \to
{\mbox{Re}}\,\O_0|_{\S^\infty},
$$
and the tangent cone $\S^\infty$ is of integer multiplicity by the 
integral compactness theorem of Allard ([A] and [Si1] 42.8). It
follows that ${\mbox{Re}}\,\O_0|_{\S^\infty}>0$, which implies
that the tangent cone $\S^\infty$ is orientable. Since $\S^\infty$
is of integer multiplicity, we have that $d\mu^\infty
=\eta(x){\cal H}^n$ where $\eta(x)$ is a locally ${\cal
H}^n$-integrable positive integer-valued function. So the cone is
an integral current (see Definition 27.1 in [Si1]).

We now show that the volume density of $\S^\infty$ at $X_0$ is
greater than 1. Otherwise, we would have
$$
\lim_{\rho\to 0}\fr{1}{\o_n\rho^n}\mu^\infty(B_\rho(0))\leq 1
$$
where $\o_n$ is the volume of the unit $n$-ball in ${\bf R}^n$:
$$
\o_n=\frac{\pi^{n/2}}{\Gamma(\frac{n}{2}+1)}.
$$
It then follows from (\ref{smon}) that for any $\e1>0$, there are
$\d1>0$ and $k_0>0$ such that for any $0<\rho<2\d1$ and $k>k_0$,
 \be\label{rec11}
\rho^{-n}\mu_{0-r^2}^k(B_\rho(\xi))< \o_n(1+\e1)
\ee
for any fixed $r>0$. The choice of $r$ will be based on
the following observation. Set
$$
\Phi(F,\x2_0,t_0,t)=\int_{\S_t}\phi(F)\frac{1}{(4\pi(t_0-t))^{n/2}}
e^{-\frac{|F-\x2_0|^2}{4(t_0-t)}}d\mu_t
$$
where $\phi$ is supported in $B_\delta(0)$ and $0\leq\phi\leq 1,\phi\equiv 1$ in
$B_{\delta/2}(0)$.
Then we have
\bs
\Phi(F_k,0,0,0-r^2)&\leq &\fr{1}{(4\pi r^2)^{n/2}}\int_0^{\d1}\int_{\part
B_\rho(0)\cap\g1_{0-r^2}^k}e^{-\fr{\rho^2}{4r^2}}d\mu_{0-r^2}^k\\
&\leq &\fr{1}{(4\pi r^2)^{n/2}}\int_0^{\d1} e^{-\fr{\rho^2}{4r^2}}\int_{
\partial B_\rho(0)\cap\g1_{0-r^2}^k}d\mu_{0-r^2}^k\\
&\leq &\fr{1}{(4\pi r^2)^{n/2}}\int_0^{\d1}
e^{-\fr{\rho^2}{4r^2}}\fr{\rho}{2r^2}{\rm Vol}(
B_\rho(0)\cap\g1_{0-r^2}^k)d\rho\\
&&+\frac{1}{(4\pi r^2)^{n/2}}e^{-\frac{\delta^2}{4r62}}
\mbox{Vol}(B_\rho(0)\cap\S^k_{0-r^2})\\
&\leq &\fr{(1+\epsilon)\o_n}{(4\pi r^2)^{n/2}}\int_0^\d1
e^{-\fr{\rho^2}{4r^2}}\fr{\rho^{n+1}}{2r^2}d\rho + o(r)~~{\rm
by}~~(\ref{rec11})~~{\rm and}~~(\ref{fm})\\
&=&\frac{(1+\epsilon)\o_n}{\pi^{n/2}}\int^{\frac{\delta^2}{4r^2}}_0
e^{-x}x^{\frac{n}{2}}dx +o(r)\\
&\leq& 1+\epsilon+o(r)
\es
because $\Gamma(\frac{n}{2}+1)=\int^\infty_0e^{-x}x^{\frac{n}{2}}dx$.
Choosing $r>0$ sufficiently small, we therefore have
$$
\Phi (F,X_0,T,T-\lmd_k^{-2}r^2)=\Phi(F_k,0,0,0-r^2)\leq 1+\e1.
$$
Now by White's local regularity theorem ([Wh1] Theorem 3.1
and Theorem 4.1, also see [E2]), $(X_0,T)$ could not be
a singular point of the mean curvature flow. This is a contradiction.
\hfill Q.E.D.

\section{Flatness of $\lmd$-cone in dimension 2}

Regularity of the $\lambda$ tangent cone can be greatly improved in the
2-dimensional case: $\dim_{\bf C}M=2$.

\begin{theorem}
Let $(M,\O)$ be a compact Calabi-Yau manifold and let $\S_0$ be a
compact Lagrangian surface in $M$ which is almost calibrated by
${\mbox{Re}}\O$. If $0<T<\infty$ is the first blow-up time of a
mean curvature flow of $\S_0$ in $M$, then the $\lmd$ tangent cone
at $(X_0,T)$ consists of a finite union (but more than two)
2-planes in ${\bf R}^4$ which are complex in a complex structure
on ${\bf R}^4$.
\end{theorem}
{\it Proof:} We use the same notation as that in the proof of
Theorem \ref{stat}, we shall show that $\t1_0$ is constant ${\cal
H}^2$ a.e. on $\g1^\infty$. To do so, we claim that for any $r>0$,
$\xi_1,\xi_2\in \g1_t^k\cap B_{R/2}(0)$ the following holds
$$\left|\fr{1}{{\rm Vol}(B_r(\xi_1)\cap\g1_t^k)}
\int_{B_r(\xi_1)\cap\g1_t^k}\cos\theta_k d\mu_t^k
-\fr{1}{{\rm Vol}(B_r(\xi_2)\cap\g1_t^k)}
\int_{B_r(\xi_2)\cap\g1_t^k}\cos\theta_k d\mu_t^k \right|
$$
\begin{equation}\label{t1c}
\leq \fr{C_1(r)}{{\rm Vol}(B_r(\xi_1)\cap\g1_t^k)}\cdot
\fr{C_2(r)}{{\rm Vol}(B_r(\xi_2)\cap\g1_t^k)}\int_{B_R(0)\cap\g1_t^k}
\left|\btd\cos\theta_k\right|d\mu_t^k,
\end{equation}
where $B_r(\xi_i)$, $i=1,2$, are the 4-dimensional balls in $M$.
To prove (\ref{t1c}), let us first recall the isoperimetric
inequality on $\g1_t^k$ (c.f. [HSp] and [MS]): let $B^k_\rho(p)$
be the geodesic ball in $\Sigma^k_t$, with radius $\rho$ and
center $p$, then \bs {\rm Vol}(B^k_\rho(p))&\leq & C\left({\rm
length }(\partial
(B^k_\rho(p)))+\int_{B^k_\rho(p)}|\h2_k|d\mu_t^k\right)^2\\
&\leq & C\left({\rm length }(\partial
(B^k_\rho(p)))+\left(\int_{B^k_\rho(p)}|\h2_k|^2d\mu_t^k
\right)^{1/2}{\rm
Vol}^{1/2}(B^k_\rho(p))\right)^2, \es for any $p\in \g1_t^k$, and
any $\rho >0$, where $C$ does not depend on $k$, $\rho$, and $p$.
By Proposition \ref{main2.1}, we have
$$
\int_{B^k_\rho(p)}|\h2_k|^2d\mu_t^k \to 0~{\rm as}~k\to\infty .
$$
So, for $k$ sufficiently large, we obtain:
$$
{\rm Vol}(B^k_\rho(p))\leq  C\left({\rm length }(\partial
(B^k_\rho(p)))\right)^2.
$$
In particular, for $k$ sufficiently large, the isoperimetric
inequality implies
\begin{equation}\label{iso}
{\rm Vol}(B^k_\rho(p))\geq C\rho^2,
\end{equation}
where $C$ is a positive constant independent of $k,\rho$ and $p$.

Suppose that the diameter of $B_r(\xi)\cap\g1_t^k$ is $d_k(\xi)$.
Then
 \bs
 Cr^2 &\geq & \int_{B_r(\xi)\cap\g1_t^k}d\mu_t^k~~~{\rm by}~(\ref{fm})\\
 &=&\int_0^{d_k(\xi)/2}\int_{\partial B^k_\rho(p)}d\sigma d\rho ~~~{\rm
 for~some}~p\in\g1_t^k\\
 &\geq & c \int_0^{d_k(\xi)/2}{\rm
 Vol}^{1/2}(B^k_\rho(p))d\rho+o(1),~o(1)\to 0~{\rm as}~k\to\infty\\
 &\geq& c\int^{d_k(\xi)/2}_0 C\rho d\rho+o(1)~~{\rm by}~(\ref{iso})\\
 &\geq & cd_k(\xi)^2+o(1).
 \es
We therefore have, for any $\xi$, \be \label{t1c2} d_k(\xi)\leq Cr+o(1) \ee
where the constant $C$ is independent of $\xi$ and $k$.

For any fixed $\eta\in B_r(\xi_2)\cap\g1_t^k$ and any $\xi\in
B_r(\xi_1)\cap\g1_t^k$, we choose a geodesic $l_{\eta\xi}$
connecting $\eta$ and $\xi$, call it a ray from $\eta$ to $\xi$.
Take an open tubular neighborhood $U(l_{\eta\xi})$ of
$l_{\eta\xi}$ in $\g1_t^k$. Within this neighborhood
$U(l_{\eta\xi})$, we call the line in the normal direction of the
ray $l_{\eta\xi}$ the normal line which we denote by
$n(l_{\eta\xi})$. It is clear that \be\label{by parts}
\cos\theta_k(\xi)-\cos\theta_k(\eta)= \int_{l_{\eta\xi}}\partial_l
\cos\theta_k dl \ee where $dl$ is the arc-length element of
$l_{\eta\xi}$. Choose $r$ small enough so that
$B_r(\xi_1)\cap\g1^k_t$ is contained in $U(l_{\eta\xi_1})$.
Keeping $\eta$ fixed and integrating (\ref{by parts}) with respect
to the variable $\xi$, first along the normal direction
$n(l_{\eta\xi_1})$ and then on the ray direction $l_{\eta\xi_1}$,
we have
\begin{eqnarray}\label{average}
\lefteqn{\left|\fr{1}{{\rm
 Vol}(B_r(\xi_1)\cap\g1_t^k)}\int_{B_r(\xi_1)\cap\g1_t^k}\cos\theta_k(\xi)
 d\mu_t^k-\cos\theta_k(\eta)\right|}\nonumber\\
 &\leq & \fr{1}{{\rm
Vol}(B_r(\xi_1)\cap\g1_t^k)}\int_{0}^{d_k(\xi_1)}\int_{n(l_{\eta\xi_1})}\int_{l_{\eta\xi}}
\left |\btd \cos\theta_k
 \right |dl d n(\xi)d\rho\nonumber\\
 &\leq & \fr{1}{{\rm
Vol}(B_r(\xi_1)\cap\g1_t^k)}\int_{0}^{d_k(\xi_1)}\int_{B_R(0)}
\left |\btd \cos\theta_k
 \right |d\mu_t^kd\rho\nonumber\\
&\leq&\fr{Cr}{{\rm Vol}(B_r(\xi_1)\cap\g1_t^k)}\int_{B_R(0)}\left
|\btd \cos\theta_k \right |d\mu_t^k,
\end{eqnarray}
here in the last step we have used (\ref{t1c2}). From
(\ref{average}), integrating with respect to $\eta$ in
$B_r(\xi_2)\cap\g1_t^k$ and dividing by ${\rm
Vol}(B_r(\xi_2)\cap\g1^k_t)$, we get the desired inequality
(\ref{t1c}).

For $i=1,2$ H\"older's inequality and (\ref{fm}) lead to
$$
\int_{B_r(\xi_i)\cap\g1^k_t}\left|\btd\cos\theta_k\right|d\mu^k_t\leq
Cr\left(\int_{B_r(\xi_i)\cap\g1^k_t}\left|\btd\cos\theta_k\right|^2d\mu^k_t\right)^{1/2}.
$$
The triangle inequality implies $B^k_r(\xi_i)\subset
B_r(\xi_i)\cap\g1^k_t$ for $i=1,2$; therefore by (\ref{iso})
$$
{\rm Vol}(B_r(\xi_i)\cap \g1^k_t)\geq {\rm Vol}(B^k_r(\xi_i))\geq
Cr^2.
$$
Now first letting $k\to\infty$ in (\ref{t1c}) and using that the
right hand side of (\ref{t1c}) tends to 0 by Proposition
\ref{main2.1}, and then letting $r\to 0$, we conclude that $\t1$
is constant ${\cal H}^2$ a.e. on $\g1^\infty$.

The $(2,0)$-form $\O_0$ is fixed by $\O(\x2_0)$ hence it has unit
length. In the complex structure $J_{\x2_0}$ on ${\bf R}^4$,
$\O_0=dz_1\wedge dz_2$.  We define
a new complex structure $J^*$ on ${\bf R}^4$:
$$ J^*(\p/\p x_1) = \theta_0 (\p/\p y_1),~~
  J^*(\p/\p y_1) = -1/\theta_0 (\p/\p x_1),
$$
$$ J^*(\p/\p x_2) = 1/\theta_0 (\p/\p y_2),~~
  J^*(\p/\p y_2) = -\t1_0(\p/\p x_2).
$$
In $J^*$, the complex coordinates are: $z^*_1 = x_1 +
\sqrt{-1}\theta_0^{-1}y_1$, $z^*_2 = \theta_0^{-1}x_2 +
\sqrt{-1}y_2$. Then $\O_0^{*}=dz^*_1\wedge dz^*_2$ satisfies that
${\mbox{Re}}\,\O_0^{*}|_{\S^\infty}= d\mu^\infty$.

We can further choose a new complex structure $J'$ on ${\bf R}^4$
such that $\O_0^*$ is of type $(1,1)$ in $J'$. In fact, if we express
$J^*$ in the local coordinates $x_1,\theta_0^{-1}y_1,\theta_0^{-1}x_2, y_2$
by
$$
J^*=\left(\begin{array}{clcr}I&0\\0&I\end{array}\right),\,\,\,
~~~{\mbox{with}}~~
I=\left(\begin{array}{clcr}0&1\\-1&0\end{array}\right),
$$
then we can take
$$
J'=\left(\begin{array}{clcr}I&0\\0&-I\end{array}\right).
$$
Therefore $\S^\infty$ is a stationary rectifiable current of type
$(1,1)$ with respect to the complex structure $J'$. By
Harvey-Shiffman's Theorem 2.1 in [HS], $\S^\infty$ is a
$J'$-holomorphic subvariety of complex dimension one. It then
follows that the singular locus ${\mathcal S}$ of $\S^\infty$
consists of isolated points.

Without loss of any generality, we may assume
$0\in\g1^\infty$ where $0$ is the origin of ${\bf R}^4$.
In fact, if not, $\g1^\infty$ would
move to infinity, then we would have
 $$
\Phi (F,X_0,T,T-\lmd_k^{-2}r^2)=\Phi(F_k,0,0,0-r^2)\to 0~{\rm
as}~k\to\infty .
$$
But White's regularity theorem then implies that $(X_0,T)$ is a
regular point. This is impossible.

There is a sequence of points $X_k\in\g1_t^k$ satisfying $X_k\to
0$ as $k\to\infty$. By Proposition \ref{main2.1}, for any $s_1$ and $s_2$ with
$-\infty<s_1<s_2<0$ and any $R>0$, we have
$$
\int_{s_1}^{s_2}\int_{\g1_t^k\cap
B_R(0)}\left|F_k^{\perp}\right|^2d\mu_t^k dt\to 0~~{\rm
as}~~k\to\infty.
$$
Thus, by (\ref{fm})
 \bs
\lefteqn{\lim_{k\to\infty}\int_{s_1}^{s_2}\int_{\g1_t^k\cap
B_R(0)}\left|(F_k-X_k)^{\perp}\right|^2d\mu_t^k dt}\\
&\leq & 2\lim_{k\to\infty}\int_{s_1}^{s_2}\int_{\g1_t^k\cap
B_R(0)}\left|F_k^{\perp}\right|^2d\mu_t^k dt+ C(s_2-s_1)R^2
\lim_{k\to\infty}|X_k|^2\\
&=& 0.
 \es

Let us denote the tangent spaces of $\g1^k_t$ at the point $F_k(x,t)$
and of $\g1^\infty$ at the point $F^\infty(x,t)$ by $T\g1_t^k$ and $T\g1^\infty$
respectively. It is clear that
$$
(F_k-X_k)^{\perp}={\rm dist}~(X_k, T{\g1_t^k}),
$$
and
$$
(F_\infty)^{\perp}={\rm dist}~(0, T{\g1^\infty}).
$$
By Allard's compactness theorem, we have
\bs \int_{s_1}^{s_2}\int_{\g1^\infty\cap
B_R(0)}\left|(F_\infty)^{\perp}\right|^2d\mu^\infty dt &=&
\int_{s_1}^{s_2}\int_{\g1^\infty\cap B_R(0)}\left|{\rm dist}~(0,
T{\g1^\infty})\right|^2d\mu^\infty dt\\
&=& \lim_{k\to\infty}\int_{s_1}^{s_2}\int_{\g1_t^k\cap
B_R(0)}\left|{\rm dist}~(X_k, T{\g1_t^k})\right|^2d\mu_t^k dt\\
&=& \lim_{k\to\infty}\int_{s_1}^{s_2}\int_{\g1_t^k\cap
B_R(0)}\left|(F_k-X_k)^{\perp}\right|^2d\mu_t^k dt\\
&=&0.
\es
Therefore $F^\perp_\infty\equiv 0$. Differentiating $\la F_\infty,v_\al\ra=0$,
inner product is taken in ${\bf R}^4$, leads to
$$
0=\la \p_i F_\infty,v_\al\ra+\la F_\infty ,\p_i v_\al\ra =\la
F_\infty ,\p_i v_\al\ra.
$$
Because $\p_i F_\infty$ is tangential to $\Sigma^\infty$,
by Weingarten's equation we observe
$$
(h_\infty)_{ij}^\al\la\f2_\infty,e_j\ra=0~~{\rm for}~~{\rm
all}~~\al,~i=1,2.
$$
Since either $\langle F_\infty, e_1\rangle\not=0$ or 
$\langle F_\infty, e_2\rangle\not=0$, we conclude $\det(h^\al_{ij})=0$. 
Recall $h^\al_{11}+h^\al_{22}=0$. It then follows
$h^\al_{ij}=0,~{\rm for}~i,j,\al=1,2$. Now we conclude that
$\S^\infty$ consists of flat 2-planes. \hfill Q.E.D.

\section{Tangent cones from a time dependent scaling}

In this section, we consider the tangent cones which arise from
the rescaled submanifold $\tilde{\g1}_s$ defined by
\begin{equation}
\tilde{F}(\cdot ,s)=\fr{1}{\sqrt{2(T-t)}}F(\cdot ,t),
\end{equation}
where $s=-\fr{1}{2}\log (T-t)$, $c_0\leq s<\infty $. Here we
choose the coordinates so that $X_0=0$. Rescaling of this type was
used by Huisken [H2] to distinguish Type I and Type II
singularities for mean curvature flows. Denote the rescaled
submanifold by $\tilde{\g1}_s$. From the evolution equation of $F$
we derive the flow equation for $\tilde{F}$
\begin{equation}\label{meaneqn2}
\fr{\p}{\p s}\tilde{F}(x,s)=\tilde{\h2}(x,s)+\tilde{\f2}(x,s).
\end{equation}
It is clear that
\begin{eqnarray*}
\cos\tilde{\theta}(x,s)&=&\cos\theta(x,t),\\
|\tilde{\h2}|^2(x,s)&=&2(T-t)|\h2|^2(x,t),\\
|\tilde{\a2}|^2(x,s)&=&2(T-t)|\a2 |^2(x,t).
\end{eqnarray*}
We set $\tilde{v}(x,t)=\cos\tilde{\theta}(x,s)$.

\begin{lemma}\label{tkae}
Assume that $(M,\O)$ is a compact Calabi-Yau manifold and $\g1_t$ evolves by
a mean curvature flow in $M$ with the initial submanifold $\S_0$ being Lagrangian
and almost calibrated by ${\mbox{Re}}\,\O$. Then
\be\label{tkae1}
\left(\fr{\p}{\p s}-\tilde{\btu} \right)\tilde{v}(x,s) = |\tilde{\h2}|^2\tilde{v}(x,s).\ee
\end{lemma}
{\it Proof:} One can check directly that
$$
\left(\fr{\p}{\p s}-\tilde{\btu}
\right)\cos\tilde{\al}(x,s)=2(T-t)\left(\fr{\p}{\p t}-\btu
\right)\cos\al(x,t).
$$
It follows that
\bs
\left(\fr{\p}{\p s}-\tilde{\btu}\right)\tilde{v}(x,s)&=&
2(T-t)\left(\fr{\p}{\p t}-\btu \right)v(x,t)\\
&\geq&2(T-t)\left|\h2\right|^2 v(x,t)\\
&=& |\tilde{\h2}|^2\tilde{v}(x,s).
\es
This proves the lemma. \hfill{Q.E.D.}

Next, we shall derive the corresponding weighted monotonicity formula
for the scaled flow. By (\ref{tkae1}), we have
$$
\left(\fr{\p}{\p s}-\tilde{\btu} \right)\fr{1}{\tilde{v}} =
-\frac{|\tilde{\h2}|^2}{\tilde{v}}-\fr{2|\tilde{\btd}{\tilde{v}}|^2}{\tilde{v}^3}.
$$
Let
$$
\tilde{\rho}(X)={\rm exp}\left(-\fr{1}{2}|X|^2\right),
$$
$$
\Psi (s)=\int_{\tilde{\g1}_s}\fr{1}{\tilde{v}}\phi\tilde{\rho}
(\tilde{\f2})d\tilde{\mu}_s.
$$
\begin{lemma} There are positive constants $c_1$ and $c_2$ which
depend on $M$, $\f2_0$ and $r$ which is the constant in the
definition of $\phi$, so that the following monotonicity formula
holds \ba \label{mon3} \fr{\p}{\p s}{\rm exp} (c_1 e^{-s})\Psi (s)
&\leq&-{\rm exp} (c_1 e^{-s}) \left(
\int_{\tilde{\g1}_s}\fr{1}{\tilde{v}}\phi\tilde{\rho}
(\tilde{\f2})
\left|\tilde{\h2}+ \tilde{F}^{\perp}\right|^2d\tilde{\mu}_s \right. \no\\
&&\left.+ \int_{\tilde{\g1}_s}\fr{1}{\tilde{v}}\phi\tilde{\rho}
(\tilde{\f2})\frac{|\tilde{\h2}|^2}{2}d\tilde{\mu}_s
+\int_{\tilde{\g1}_s}\fr{2}{\tilde{v}^3}\left|\tilde{\btd}
\tilde{v}\right|^2\phi\tilde{\rho} (\tilde{\f2})d\tilde{\mu}_s
\right)\no\\
&&+c_2{\rm exp} (c_1 e^{-s}) . \ea
\end{lemma}
{\it Proof:} Note that
\begin{eqnarray*}
&&\tilde{F}(x,s)=\fr{F(x,t)}{\sqrt{2(T-t)}},\\
&&\tilde{\h2}(x,s)=\sqrt{2(T-t)}\h2(x,t),\\
&&|\tilde{\btd}\tilde{v}|^2(x,s)=2(T-t)|\btd v|^2(x,t).
\end{eqnarray*}
By the chain rule
$$\fr{\p}{\p s}=2e^{-2s}\fr{\p}{\p t}$$
and the monotonicity inequality (\ref{mon2}) for the unscaled
submanifold, we obtain the desired inequality. \hfill{Q.E.D.}

\begin{lemma}
\label{main7.1} Let $(M,\O)$ be a compact Calabi-Yau manifold. If
the initial compact submanifold $\S_0$ is Lagrangian and almost
calibrated by $\mbox{Re}\,\O$, then there is a sequence
$s_k\to\infty$ such that, for any $R>0$, \be\label{m7.2}
\int_{\tilde{\g1}_{s_k}\cap B_R(0)} |\tilde{\btd} \cos
\tilde{\theta}|^2d\tilde{\mu}_{s_k}\to 0~~{\rm as}~~k\to\infty,
\ee \be\label{m7.3} \int_{\tilde{\g1}_{s_k}\cap
B_R(0)}|\tilde{\h2}|^2d\tilde{\mu}_{s_k}\to 0~~{\rm
as}~~k\to\infty, \ee and \be\label{m7.4}
\int_{\tilde{\g1}_{s_k}\cap
B_R(0)}|\tilde{F}^{\perp}|^2d\tilde{\mu}_{s_k}\to 0~~{\rm
as}~~k\to\infty. \ee
\end{lemma}
{\it Proof:} Integrating (\ref{mon3}), we have \bs  \infty & > &
\int_{s_0}^\infty\int_{\tilde{\g1}_s}\fr{1}{\tilde{v}}\phi\tilde{\rho}
(\tilde{\f2})
\left|\tilde{\h2}+ \tilde{F}^{\perp}\right|^2d\tilde{\mu}_s ds \\
&&+
\int_{s_0}^\infty\left(\int_{\tilde{\g1}_s}\fr{1}{\tilde{v}}\phi\tilde{\rho}
(\tilde{\f2})\frac{|\tilde{\h2}|^2}{2}d\tilde{\mu}_s
+\int_{\tilde{\mu}_s}\fr{2}{\tilde{v}^3}|\tilde{\btd}
\tilde{v}|^2\phi\tilde{\rho} (\tilde{\f2})d\tilde{\mu}_s
\right)ds. \es
Hence there is a sequence $s_k\to\infty$, such that as $k\to\infty$
$$
\int_{\tilde{\g1}_{s_k}}\fr{1}{\tilde{v}}\phi\tilde{\rho}
(\tilde{\f2})\frac{|\tilde{\h2}|^2}{2}d\tilde{\mu}_{s_k}\to 0,
$$
$$
\int_{\tilde{\g1}_{s_k}}\fr{2}{\tilde{v}^3}|\tilde{\btd}
\tilde{v}|^2\phi\tilde{\rho} (\tilde{\f2})d\tilde{\mu}_{s_k}\to 0,
$$
and
$$
\int_{\tilde{\g1}_{s_k}}\fr{1}{\tilde{v}}\phi\tilde{\rho}
(\tilde{\f2}) \left|\tilde{\h2}+
\tilde{F}^{\perp}\right|^2d\tilde{\mu}_{s_k}\to 0.
$$
Since $\tilde{v}$ has a positive lower bound,
the proposition now follows. \hfill Q.E.D.

The proof of the following lemma is essentially the same as the
one for Proposition \ref{rectif}, except there are two parameters
$\lmd,t$ for the $\lmd$ tangent cones but only one parameter $t$
for the time dependent tangent cones. Note that the alternative proof
given in [CL1] using the isoperimetric inequality only works in dimension 2.

\begin{lemma}\label{texit}
There is a subsequence of $s_k$, which we also denote by $s_k$,
such that $(\tilde{\g1}_{s_k},d\tilde{\mu}_{s_k})\to
(\tilde{\g1}_\infty,d\tilde{\mu}_\infty)$ in the sense of measures.
And $(\tilde{\g1}_\infty,d\tilde{\mu}_\infty)$ is ${\mathcal H}^n$-rectifiable.
\end{lemma}
{\it Proof:} To show the subconvergence, it suffices to show that,
for any $R>0$,
\be\label{tfm}
\tilde{\mu}_{s_k}(\tilde{\g1}_{s_k}\cap B_R(0))\leq CR^n,
\ee
where $B_R(0)$ is a metric ball in ${\bf R}^{2n}$, $C>0$ is
independent of $k$. Direct calculation leads to
\bs
\tilde{\mu}_{s_k}(\tilde{\g1}_{s_k}\cap B_R(0))
&=&(2(T-t))^{-n/2}\int_{\g1_{T-e^{2s_k}}\cap B_{\sqrt{2(T-t)}R}(0)}d\mu_t\\
&=& R^n\left(\sqrt{2}e^{-s_k}R\right)^{-n}\int_{\g1_{T-e^{2s_k}}\cap
B_{\sqrt{2}e^{-s_k}R}(0)}d\mu_t\\
&\leq& CR^n\int_{\g1_{T-e^{2s_k}}\cap B_{\sqrt{2(T-t)}R}(0)}\frac{1}{v}
\frac{1}{(4\pi)^{n/2}(\sqrt{2}e^{-s_k}R)^n}e^{-\frac{|X-X_0|^2}{4\sqrt{2}e^{-s_k}R}}d\mu_t\\
&\leq & CR^n\Psi\left(0,T+(\sqrt{2}e^{-s_k}R)^{2}-e^{2s_k},
T-e^{2s_k}\right)
\es
By the monotonicity inequality (\ref{mon2}), we have
\bs
\tilde{\mu}_{s_k}(\tilde{\g1}_{s_k}\cap B_R(0))
&\leq&CR^n\Phi(0,T+(\sqrt{2}e^{-s_k}R)^{2}-e^{2s_k}, T/2)+CR^n\\
&\leq& \frac{C\mu_{T/2}(\g1_{T/2})}{T^{n/2}\min_{\S_0}v}R^n+CR^n.
\es
Since volume is non-increasing along mean curvature flow, we see
$$
\tilde{\mu}_{s_k}(\tilde{\g1}_{s_k}\cap B_R(0))\leq CR^n.
$$

We now prove that $(\tilde{\g1}_\infty,d\tilde{\mu}_\infty)$ is
${\mathcal H}^n$-rectifiable. For any $\xi\in\tilde{\g1}_\infty$,
choose $\xi_k\in\tilde{\g1}_{s_k}$ with $\xi_k\to\xi$ as
$k\to\infty$. By the monotonicity identity (17.4) in [Si1], we
have \ba\label{tsmon}
\sigma^{-n}\tilde{\mu}_{s_k}(B_\sigma(\xi_k))&=&
\rho^{-n}\tilde{\mu}_{s_k}(B_\rho(\xi_k))
-\int_{B_\rho(\xi_k)\setminus B_\sigma(\xi_k)}\fr{|D^\perp
r|^2}{r^n} d\tilde{\mu}_{s_k}
\no\\
&& -\fr{1}{n}\int_{B_\rho(\xi_k)}(x-\xi_k)\cdot \tilde{{\bf H}_k}
\left(\fr{1}{r_\sigma^n}-\fr{1}{\rho^n}\right)d\tilde{\mu}_{s_k},
\ea for all $0<\sigma \leq \rho$, where
$\tilde{\mu}_{s_k}(B_\sigma(\xi_k))$ is the area of
$\tilde{\g1}_{s_k}\cap B_\sigma(\xi_k)$, $r_\sigma =\max
\{r,\sigma\}$ and $D^\perp r$ denotes the orthogonal projection of
$Dr$ (which is a vector of length 1) onto
$(T_{\xi_k}\tilde{\g1}_{s_k})^\perp$. Letting $k\to\infty$, by
Lemma \ref{main7.1}, we have
$$
\sigma^{-n}\tilde{\mu}_\infty(B_\sigma(\xi))\leq
\rho^{-n}\tilde{\mu}_\infty(B_\rho(\xi)),
$$
for all $0<\sigma \leq \rho$. Therefore, $\lim_{\rho\to
0}\rho^{-n}\tilde{\mu}_\infty(B_\rho(\xi))$ exists and is finite by (\ref{tfm}).

By converting $s$ to $t$, the argument for the positive lower bound
of the volume density in the proof of
Proposition \ref{rectif} carries over to the present situation.

We conclude that $\lim_{\rho\to
0}\rho^{-n}\tilde{\mu}_\infty(B_\rho(\xi))$ exists and for
${\mathcal H}^n$ almost all $\xi\in \tilde{\g1}_\infty$,
\begin{equation}\label{tden}
0<C\leq\lim_{\rho\to
0}\rho^{-n}\tilde{\mu}_\infty(B_\rho(\xi))<\infty.
\end{equation}
Priess's theorem in [P] then asserts the ${\mathcal H}^n$-rectifiability of
$(\tilde{\g1}_\infty,d\tilde{\mu}_\infty)$. \hfill Q.E.D.

\begin{definition}
\label{bbtt} {\em  We call
$(\tilde{\g1}_\infty,d\tilde{\mu}_\infty)$ obtained in
Lemma \ref{texit} {\it a tangent cone  of the mean curvature flow $\g1_t$ at
$(X_0,T)$ in the time dependent scaling.}}
\end{definition}

With the lemmas established in this section, by using arguments
completely similar to those for the $\lmd$ tangent cones in the previous sections,
we can prove

\begin{theorem}\label{tflat}
Let $(M,\O)$ be a compact Calabi-Yau manifold. If the initial
compact submanifold $\S_0$ is Lagrangian and almost calibrated by
${\mbox{Re}}\,\O$ and $T>0$ is the first blow-up time of the mean
curvature flow, then the tangent cone $\tilde{\g1}_\infty$ of the
mean curvature flow at $(X_0,T)$ coming from time dependent
scaling is a rectifiable stationary Lagrangian current with integer multiplicity 
in ${\bf R}^{2n}$.
Moreover, if $M$ is of complex 2-dimensional, then
$\tilde{\S}_\infty$ consists of a finitely many (more than 1)
2-planes in ${\bf R}^4$ which are complex in a complex structure
on ${\bf R}^4$.
\end{theorem}

The result below can also be found in [Wa].
\begin{corollary}
If the initial compact submanifold $\S_0$ is Lagrangian and is almost calibrated in
a compact Calabi-Yau manifold $(M,\O)$, then mean curvature flow does not develop
Type I singularity.
\end{corollary}
{\it Proof:} Let $X_0$ be a Type I singularity at $T<\infty$ and set $\lmd=\max_{\S_t}|\a2|^2$.
The $\lmd$ tangent cone $\S_\infty$ is smooth if $T$ is a Type I singularity.
Therefore $\S_\infty$ is a smooth minimal Lagrangian submanifold in ${\bf C}^n$ by Theorem \ref{tflat}.
Because $\S_\infty$ is smooth, (\ref{m2.4}) implies $F_\infty^\perp\equiv 0$ everywhere. 
The monotonicity identity (\ref{smon}) then implies $\sigma^{-n}\mu(\S_\infty\cap B_\sigma(0))$ 
is a constant independent of $\sigma$, and the volume density ratio at 0 is one due to 
the smoothness of $\S_\infty$, 
so $\S_\infty$ is a flat linear subspace of ${\bf R}^{2n}$. But the second
fundamental form of $\S_\infty$ has length one at $0$ according to the blow-up process,
and the contradiction rules out any Type I singularity. \hfill{Q.E.D.}

\vspace{.2in}
\begin{center}
{\large\bf REFERENCES}
\end{center}
\footnotesize
\begin{description}
\item[{[A]}] {W. Allard, First variation of a varifold, Annals
of Math. {\bf 95} (1972), 419-491.}
\item[{[B]}] {K. Brakke, The motion of a surface by its mean curvature,
Princeton Univ. Press, 1978.}
\item[{[CL1]}] {J. Chen and J. Li, Singularities of codimension two
mean curvature flow of symplectic surfaces, preprint (2002).}
\item[{[CL2]}] {J. Chen and J. Li, Mean curvature flow of surface in 4-manifolds,
Adv. Math. {\bf 163} (2001), 287-309.}
\item[{[CLT]}] {J. Chen, J. Li and G. Tian, Two-dimensional graphs
moving by mean curvature flow, Acta Math. Sinica, English Series
Vol.{\bf 18}, No.2 (2002), 209-224.}
\item[{[CT]}] {J. Chen and G. Tian, Moving symplectic curves in
K\"ahler-Einstein surfaces, Acta Math. Sinica, English Series,
{\bf 16} (4), (2000), 541-548.}
\item[{[E1]}] {K. Ecker, On regularity for mean curvature flow of hypersurfaces,
Calc. Var. 3(1995), 107-126.}
\item[{[E2]}] {K. Ecker, A local monotonicity formula for mean curvature flow, Ann. Math.,
154(2001), 503-525.}
\item[{[HL]}] {R. Harvey and H.B. Lawson, Calibrated geometries,
Acta Math. {\bf 148} (1982), 47-157.}
\item[{[HS]}] {R. Harvey and B. Shiffman, A characterization of
holomorphic chains, Acta Math. {\bf 148} (1982), 47--157.}
\item[{[HSp]}] {D. Hoffman and J. Spruck, Sobolev and isoperimetric
inequalities for Riemannian submanifolds, Comm. Pure Appl. Math.
{\bf 27} (1974), 715-727.}
\item[{[H1]}] {G. Huisken, Asymptotic behavior for singularities
of the mean curvature flow, J. Diff. Geom. {\bf 31} (1990),
285-299.}
\item[{[H2]}] {G. Huisken, Flow by mean curvature of convex surfaces into
spheres, J. Diff. Geom. {\bf 20} (1984), 237-266.}
\item[{[H3]}] {G. Huisken, Contracting convex hypersurfaces in Riemannian
manifolds by their mean curvature, Invent. math. {\bf 84} (1986),
463-480.}
\item[{[HS1]}] {G. Huisken and C. Sinestrari, Convexity
estimates for mean curvature flow and singularities of mean convex
surfaces. Acta Math. {\bf 183} (1999), no. 1, 45-70.}
\item[{[HS2]}] {G. Huisken and C. Sinestrari, Mean curvature flow
singularities for mean convex surfaces. Calc. Var. Partial
Differential Equations {\bf 8} (1999), no. 1, 1-14.}
\item[{[I1]}] {T. Ilmanen, Singularity of mean curvature flow of surfaces,
preprint.}
\item[{[I2]}] {T. Ilmanen, Elliptic Regularization and Partial Regularity
for Motion by Mean Curvature, Memoirs of the Amer. Math. Soc.,
520, 1994.}
\item[{[MS]}] {J.H. Michael and L. Simon, Sobolev and mean-valued
inequalities on generalized submanifolds of $R^n$, Comm. Pure
Appl. Math. {\bf 26} (1973), 361-379.}
\item[{[P]}] {D. Priess, Geometry of measures in $R^n$; Distribution,
rectifiability, and densities; Ann. of Math., {\bf 125} (1987),
537-643.}
\item[{[ScW]}] {R. Schoen and J. Wolfson, Minimizing area among Lagrangian
surfaces: the mapping problem, J. Differential Geom. {\bf 58} (2001), 1-86.}
\item[{[Si1]}] {L. Simon, Lectures on Geometric Measure Theory,
Proc. Center Math. Anal. 3 (1983), Australian National Univ.
Press.}
\item[{[Sm1]}] {K. Smoczyk, Der Lagrangesche mittlere
Kruemmungsfluss. Univ. Leipzig (Habil.-Schr.), 102 S. 2000.}
\item[{[Sm2]}] {K. Smoczyk, Harnack inequality for the Lagrangian
mean curvature flow. Calc. Var. PDE, {\bf 8} (1999), 247-258.}
\item[{[Sm3]}] {K. Smoczyk, Angle theorems for the Lagrangian mean
curvature flow, preprint 2001.}
\item[{[TY]}] {R. Thomas and S.T. Yau, Special Lagrangians,
stable bundles and mean curvature flow, math.DG/0104197 (2001).}
\item[{[Wa]}] {M.-T. Wang, Mean curvature flow of surfaces in
Einstein four manifolds, J. Diff. Geom. {\bf 57} (2001), 301-338.}
\item[{[Wh1]}] {B. White, A local regularity theorem for classical
mean curvature flow, preprint (2000).}
\item[{[Wh2]}] {B. White, Stratification of minimal surfaces, mean curvature flows, and
harmonic maps, J. Reine Angew. Math. {\bf 488} (1997), 1-35.}
\item[{[Wh3]}] {B. White, The size of the singular set in mean curvature flow
of mean-convex sets, J. Amer. Math. Soc. {\bf 13} (2000), 665-695.}
\end{description}

\end{document}